\documentclass[12pt,a4paper]{article}
\usepackage{indentfirst}
\usepackage{amsmath,amssymb,amsfonts,amsthm,graphics}
\usepackage{makeidx}
\newtheorem{thm}{Theorem}[section]

\theoremstyle{definition}
\newtheorem{defn}{Definition}[section]

\def\-{\mbox{--}}

\newtheorem{pro}{Proposition}[section]
\newtheorem{lem}{Lemma}[section]

\begin{document}
\title{\Large\bf Rainbow vertex-connection number\\
of 2-connected graphs\footnote{Supported by NSFC No.11071130.}}
\author{\small  Xueliang~Li, Sujuan~Liu\\
\small Center for Combinatorics and LPMC-TJKLC\\
\small Nankai University, Tianjin 300071, China\\
\small lxl@nankai.edu.cn; sjliu0529@126.com}
\date{}
\maketitle
\begin{abstract}
The {\em rainbow vertex-connection number}, $rvc(G)$, of a connected
graph $G$ is the minimum number of colors needed to color its
vertices such that every pair of vertices is connected by at least
one path whose internal vertices have distinct colors. In this paper
we first determine the rainbow vertex-connection number of cycle
$C_n$ of order $n\geq 3$, and then, based on it, prove that for any
2-connected graph $G$, $rvc(G)\leq rvc(C_n)$, giving a tight upper
bound for the rainbow vertex-connection. As a consequence, we show
that for a connected graph $G$ with a block decomposition $B_1, B_2,
\cdots, B_k$ and $t$ cut vertices, $rvc(G)\leq
rvc(B_1)+rvc(B_2)+\cdots +rvc(B_k)+t$.



\end{abstract}

\section{Introduction}

All graphs in this paper are finite, undirected and simple. We follow
the terminology and notation of Bondy and Murty \cite{bondy2008graph}.

Let $G$ be a vertex-colored connected graph. A path of $G$ is a
$rainbow$ $path$ if its internal vertices have distinct colors. The
vertex-colored graph $G$ is called $rainbow \ vertex$-$connected$ if
any two vertices are connected by at least one rainbow path. The
{\em rainbow vertex-connection number} of a connected graph $G$,
denoted by $rvc(G)$, is the smallest number of colors that are
needed in order to make $G$ rainbow vertex-connected. If $F$ is a
subgraph of a vertex-colored graph, $colours(F)$ denotes the set of
all colors appeared in $F$. If a rainbow vertex-coloring of $G$ uses
$k$ colors, we call it a $k$-rainbow vertex-coloring.

Let $F$ be a subgraph of a graph $G$. An $ear$ of $F$ in $G$ is a nontrivial
path whose two ends are in $F$ but whose internal vertices are not. A
nested sequence of graphs is a sequence ($G_0, G_1, \cdots, G_k$) of
graphs such that $G_i\subset G_{i+1}, 0\leq i\leq k-1$. An $ear$
$decomposition$ of a 2-connected graph $G$ is a nested sequence ($G_0,
G_1, \cdots, G_k$) of 2-connected subgraphs of $G$ satisfying the
conditions: $(1)$ $G_0$ is a cycle; $(2)$ $G_i=G_{i-1}\bigcup P_i$,
where $P_i$ is an ear of $G_{i-1}$ in $G, 1\leq i\leq k$; $(3)$
$G_k=G$. An ear with odd (resp. even) length is called an odd (resp.
even) ear.

A maximal connected subgraph without any cut vertex is called a {\em
block}. Thus, every block of a nontrivial connected graph is either
a maximal 2-connected subgraph or $K_2$. All the blocks of a graph
$G$ form a {\em block decomposition} of $G$.

Given two walks $W_1=u_0,u_1,\cdots,u_k$ and $W_2=v_0,v_1,\cdots,v_{\ell}$
such that $u_k=v_0$, we can {\em concatenate} $W_1$ and $W_2$ to get a
long walk, $W=W_1W_2=u_0,u_1,\cdots,u_k(=v_0),v_1,v_2,\cdots,v_{\ell}$.
We denote the order of a graph by $|G|$.

The concept of rainbow vertex-connection number was introduced by
Krivelevich and Yuster in \cite{krivelevich2010rainbow}. Some easy
observations about rainbow vertex-connection number include that if
$G$ is a connected graph of order $n$, then $diam(G)-1\leq
rvc(G)\leq n-2$; $rvc(G)=0$ if and only if $G$ is a complete graph
and $rvc(G)=1$ if and only if $diam(G)=2$. Krivelevich and Yuster
\cite{krivelevich2010rainbow} showed that if a connected graph $G$
has $n$ vertices and minimum degree $\delta$, then
$rvc(G)\leq11n/\delta$. In \cite{Li2010vertex}, Li and Shi improved
the bound. In \cite{Chen2010vertex}, Chen, Li and Shi studied the
computational complexity of rainbow vertex-connection and proved
that computing $rvc(G)$ is NP-hard.

In this paper the rainbow vertex-connection $rvc(C_n)$ of a cycle
$C_n(n\geq3)$ is determined. Based on it, we then prove that for any
2-connected graph $G$ of order $n(\geq3)$, $rvc(G)\leq rvc(C_n)$
which gives a tight upper bound for the rainbow vertex-connection of
a 2-connected graph. As a consequence, we show that for a connected
$G$ with a block decomposition $B_1, B_2, \cdots , B_k$ and $t$ cut
vertices, $rvc(G)\leq rvc(B_1)+ rvc(B_2)+\cdots +rvc(B_k)+t$. The
proof is constructive and hence we can give a method to construct a
rainbow vertex-coloring of the connected graph $G$ using at most
$\lceil\frac{|B_1|}{2}\rceil+\cdots+ \lceil\frac{|B_k|}{2}\rceil+t$
colors.

\section{Main results}

First of all, we need to introduce the concept of revised rainbow
vertex-coloring.

\begin{defn}
Let $G$ be a connected graph with a vertex-coloring $c$. A path $P$
of $G$ is called a {\em revised rainbow path} if all vertices of $P$
have distinct colors, or all but the end vertices of $P$ have
distinct colors, which means that only the two end-vertices of $P$
may have the same color. The vertex-coloring $c$ of $G$ is called a
{\em revised rainbow vertex-coloring} if any two vertices of $G$ are
connected by at least one revised rainbow path. The {\em revised
rainbow vertex-connection number} of a connected graph $G$, denoted
by $rvc^*(G)$, is the smallest number of colors that are needed in
order to make $G$ revised rainbow vertex-connected.
\end{defn}
Since a revised rainbow path is also a rainbow path, $rvc(G)\leq
rvc^*(G)$. At firs, we consider the rainbow vertex-connection of a
cycle.

\begin{thm}
\label{thm:cycle} Let $C_n$ be a cycle of order $n(n\geq 3)$. Then
$$
rvc(C_n)=\left\{
\begin{array}{ll}
0 & $if $ n=3; \\
1 & $if $ n=4, 5; \\
3 & $if $ n=9; \\
\lceil\frac{n}{2}\rceil-1 & $if $ n=6, 7, 8, 10, 11, 12, 13, $or $ 15; \\
\lceil\frac{n}{2}\rceil & $if $ n\geq 16 $or $ n=14.
\end{array}
\right.
$$
\end{thm}

\begin{proof}
Assume that $C_n=v_1,v_2,\cdots,v_n,v_{n+1}(=v_1)(n\geq3)$. It is obvious that
$rvc(C_3)=0$. Since $rvc(G)=1$ if and only if $diam(G)=2$, we have
$rvc(C_4)=rvc(C_5)=1$.
\begin{figure}[h,t,b,p]
\begin{center}
\scalebox{0.6}[0.6]{\includegraphics{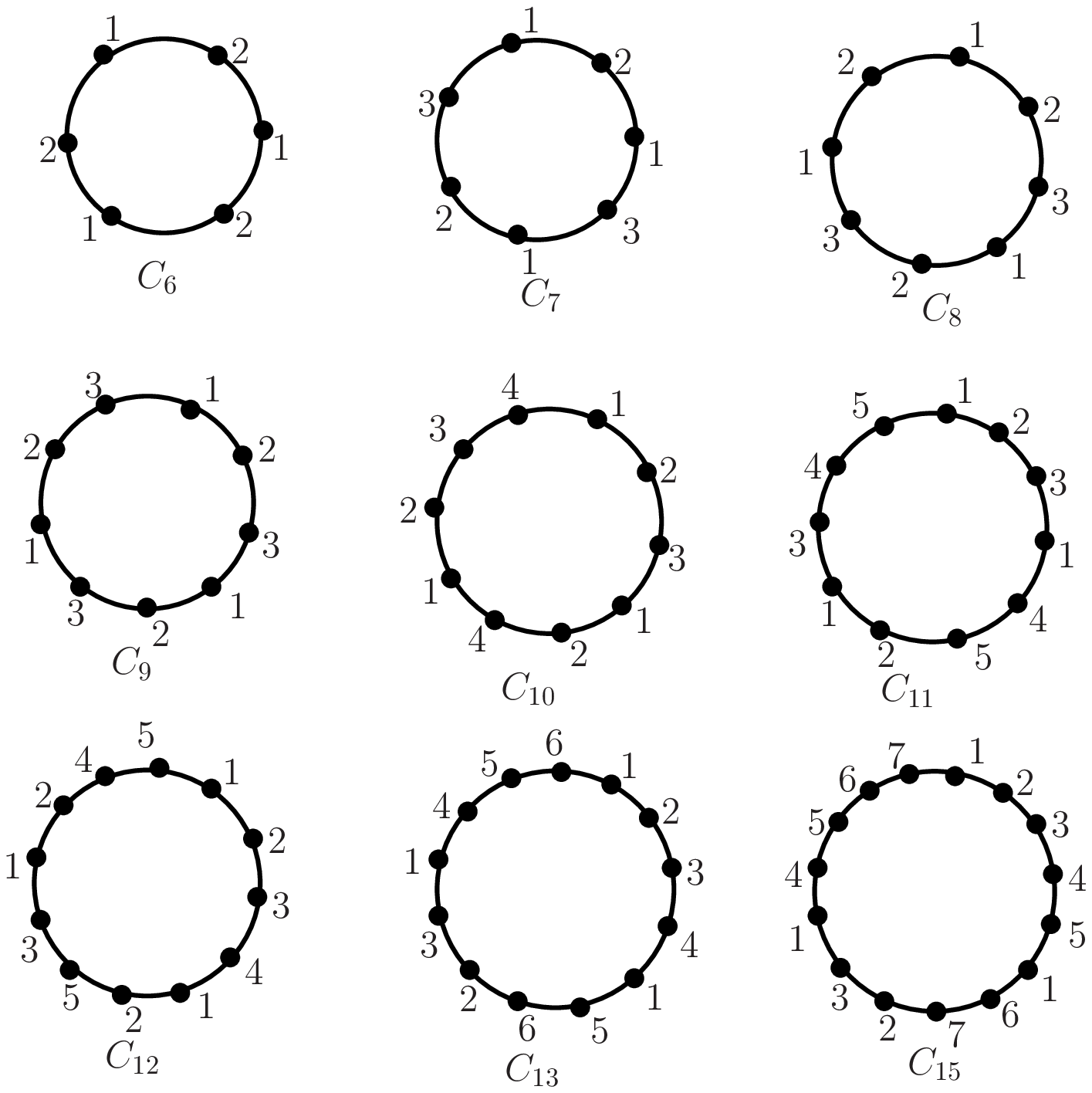}}

Figure 1. Rainbow vertex-colorings for small cycles.
\end{center}
\end{figure}

It is easy to check that the vertex-colorings of $C_n$ shown in
Fig.1 are rainbow vertex-colorings. So $rvc(C_n) \leq
\lceil\frac{n}{2}\rceil-1$ for $6 \leq n\leq 13$ and $n=15$ and
$rvc(C_9)\leq3$. Since $rvc(C_n)\geq diam(C_n)-1$, $rvc(C_n)=
\lceil\frac{n}{2}\rceil-1$ for $n=6, 8, 10, 12$ and $rvc(C_9)=3$.
For any vertex-coloring $c$ of $C_7$ using 2 colors, there exist two
adjacent vertices (say $v_1, v_2$) having the same color. Then the
two paths $P_1=v_7,v_1,v_2,v_3$ and $P_2=v_7,v_6,v_5,v_4,v_3$ on
$C_n$ are not rainbow paths, i.e., there is no rainbow path between
$v_7$ and $v_3$. So $c$ is not a rainbow vertex-coloring. Hence,
$rvc(C_7)=3$. Assume, to the contrary, that $C_n(n=11, 13, 15)$ has
a $(\lceil\frac n{2}\rceil-2)$-rainbow vertex-coloring $c$. Then
some three vertices (say $v_1, v_i, v_j\in V(C_n), 1<i<j\leq n$)
have the same color and one pair of vertices among them (say $v_1,
v_i$) has distance no more than $\lfloor\frac n{3}\rfloor$, i.e.,
$d_{C_n}(v_1, v_i)\leq\lfloor\frac n{3}\rfloor$. Let
$P=v_1,v_2,\cdots,v_i$ be the path from $v_1$ to $v_i$ with length
$d_{C_n}(v_1, v_i)$. Since $c(v_1)=c(v_i)$,
$P_1=v_n,v_1,\cdots,v_i,v_{i+1}$ is not a rainbow path. So
$P_2=C_n-P=v_n,v_{n-1},\cdots,v_{i+1}$ is the rainbow path from
$v_n$ to $v_{i+1}$. Since $\ell(C_n-P)=n-(\ell(P)+2)\geq
n-\lfloor\frac n{3}\rfloor-2 \geq\lceil\frac n{2}\rceil$ for $n=11,
13, 15$, $C_n-P$ has at least $\lceil\frac n{2}\rceil-1$ internal
vertices. So $C_n-P$ is not a rainbow path, a contradiction. Hence,
$rvc(C_n)=\lceil\frac n{2}\rceil-1$ for $n=11, 13, 15$.

In the following, we consider the rainbow vertex-connection of $C_n$ for
$n\geq 16$ or $n=14$. Suppose $C_n=v_1v_2\cdots v_nv_{n+1}(=v_1)$. Define
a vertex-coloring $c$ of $C_n$ by $c(v_i)=i$ for $1\leq i\leq
\lceil\frac{n}{2}\rceil$ and $c(v_i)=i-\lceil\frac{n}{2}\rceil$ if
$\lceil\frac{n}{2}\rceil +1\leq i \leq n$. Since for any two vertices $u, v$
of $C_n$, the path on $C_n$ with length $d_{C_n}(u, v)$ is a rainbow path,
we have $rvc(C_n)\leq
\lceil\frac{n}{2}\rceil$ for $n\geq 16$ or $n=14$.

Next we show that $rvc(C_n)\geq \lceil\frac{n}{2}\rceil$ for $n\geq
16$ or $n=14$. Assume, to the contrary, that $rvc(C_n)\leq
\lceil\frac{n}{2}\rceil-1$. Then there exists a
$(\lceil\frac{n}{2}\rceil-1)$-rainbow vertex-coloring $c_0$ of
$C_n$. Obviously, there are three vertices (say $v_1, v_i, v_j,
1\leq i< j \leq n$) of $C_n$ having the same color. And one pair of
vertices among $\{v_1, v_i, v_j\}$ (say $v_1, v_i$) satisfies that
$d_{C_n}(v_1, v_i)\leq \lfloor\frac{n}{3}\rfloor$. Without loss of
generality, assume that $P=v_1,v_2,\cdots,v_i$ is the path on $C_n$
with length $d_{C_n}(v_1, v_i)$. Now consider the vertices $v_n$ and
$v_{i+1}$. Since $v_1$ and $v_i$ have the same color, the rainbow
path between $v_n$ and $v_{i+1}$ on $C_n$ must be $C_n-P$. Since
$\ell(C_n-P)=n-(\ell(P)+2)\geq n-\lfloor\frac{n}{3}\rfloor-2$, the
number of internal vertices of $C_n-P$ is at least
$n-\lfloor\frac{n}{3}\rfloor-3$. For $n\geq 16$ or $n=14$,
$n-\lfloor\frac{n}{3}\rfloor-3 > \lceil\frac{n}{2}\rceil-1$ which
contradicts that $C_n-P$ is a rainbow path. Hence $rvc(C_n)\geq
\lceil\frac{n}{2}\rceil$ for $n\geq 16$ or $n=14$. Therefore,
$rvc(C_n)= \lceil\frac{n}{2}\rceil$ for $n\geq 16$ or $n=14$.
\end{proof}

Before proceeding, we introduce another notion.
\begin{defn}
\label{balanced}(Balanced coloring)
Let $H$ be a connected graph and $P$ an ear of $H$ such that $V(H)\bigcap
V(P)=\{a, b\}$. Assume that $P=v_1(=a),v_2,\cdots,v_s(=b)(s\geq6)$ and
$H$ has a vertex-coloring $c^\prime$. Let $x$ be a color of $c^\prime$
and $x_1, x_2,\cdots, x_{\lceil\frac s{2}\rceil-1}$ new colors.
Now define a vertex-coloring $c$ of $G:=H\bigcup P$ from $c^\prime$.
We distinguish the following cases according to the parities of
$|H|$ and $\ell(P)$.

Case 1. $|H|$ is even and $\ell(P)(=s-1)$ is odd.

$c(v)=c^\prime(v)$ for $v\in V(H)\setminus\{a,b\}$; the last $\frac s{2}$
vertices of $P$, i.e., $v_{\frac s{2}+1},\cdots,v_s$ are colored by
$c^\prime(b), x_1, \cdots,x_{\frac s{2}-1}$ in order, and if $c^\prime(a)
\neq c^\prime(b)$, then the first $\frac s{2}$ vertices of $P$, i.e.,
$v_1,\cdots,v_{\frac s{2}}$ are colored by $x_1,\cdots,x_{\frac s{2}-1},
c^\prime(a)$ in order; otherwise, by $c^\prime(a),x_1,\cdots,x_{\frac s{2}-1}$
in order.

Case 2. $|H|$ is odd and $\ell(P)$ is even.

$c(v)=c^\prime(v)$ for $v\in V(H)\setminus\{a,b\}$; the last $\lceil\frac s{2}\rceil-1$
vertices of $P$, i.e., $v_{\lceil\frac s{2}\rceil+1},\cdots,v_s$ are colored by
$c^\prime(b), x_1, \cdots,x_{\lceil\frac s{2}\rceil-2}$ in order, and if $c^\prime(a)
\neq c^\prime(b)$, then the first $\lceil\frac s{2}\rceil$ vertices of $P$, i.e.,
$v_1,\cdots,v_{\lceil\frac s{2}\rceil}$ are colored by $x_1,\cdots,x_{\lceil\frac s{2}\rceil-2},
c^\prime(a), x$ in order; otherwise, by $c^\prime(a), x_1,\cdots,$ $x_{\lceil\frac s{2}\rceil-2}$,
$x$ in order.

Case 3. $|H|$ and $\ell(P)$ are even.

First, color one vertex of $P$ by $x_{\lceil\frac s{2}\rceil-1}$.
Second, $c(v)=c^\prime(v)$ for $v\in V(H)\setminus\{a,b\}$; the last $\lceil\frac s{2}\rceil-1$
uncolored vertices of $P$ are colored by
$c^\prime(b), x_1, \cdots,x_{\lceil\frac s{2}\rceil-2}$ in order, and if $c^\prime(a)
\neq c^\prime(b)$, then the first $\lceil\frac s{2}\rceil-1$ uncolored vertices of $P$
are colored by $x_1,\cdots,x_{\lceil\frac s{2}\rceil-2},
c^\prime(a)$ in order; otherwise, by $c^\prime(a), x_1,\cdots,$ $x_{\lceil\frac s{2}\rceil-2}$
in order.

Case 4. $|H|$ and $\ell(P)$ are odd.

First, color one vertex of $P$ by $x_{\frac s{2}-1}$.
Second, $c(v)=c^\prime(v)$ for $v\in V(H)\setminus\{a,b\}$; the last $\frac s{2}-1$
uncolored vertices of $P$ are colored by
$c^\prime(b), x_1, \cdots,x_{\frac s{2}-2}$ in order, and if $c^\prime(a)
\neq c^\prime(b)$, then the first $\frac s{2}$ uncolored vertices of $P$
are colored by $x_1,\cdots,x_{\frac s{2}-2},
c^\prime(a), x$ in order; otherwise, by $c^\prime(a), x_1,\cdots,$ $x_{\frac s{2}-2}, x$
in order.

The obtained vertex-coloring $c$ of $G$ from $c^\prime$ is called a {\em balanced
coloring}.      $\Box$
\end{defn}

It is obvious that when $|G|$ is even, the balanced coloring $c$ is unique
with respect to $c^\prime$, and when $|G|$ is odd, it is not unique.

Let $G$ be a connected graph and $v$ a vertex of $G$. Assume that
$c$ is a revised rainbow vertex-coloring of $G$ and $x$ is a color
of $c$ such that $c(v)\neq x$. Define {\em a property} ($*$) of $c$
with respect to $v$ and $x$ as follows: for any vertex $u$ of $G$
which $c(u)\neq x$, there exists a revised rainbow path $P$ from $v$
to $u$ such that $x\notin colours(P)$.

\begin{pro}\label{pro1}
Let $H$ be a connected graph and $P=v_1(=a),v_2,\cdots
,v_s(=b)(s\geq6)$ an ear of $H$ such that $V(H)\bigcap V(P)=\{a,
b\}$. $H$ has a revised $\lceil\frac{|H|}{2}\rceil$-rainbow
vertex-coloring $c^\prime$ each of whose colors appears at most
twice. Moreover, when $|H|$ is odd, $c^\prime$ satisfies the
property ($*$) with respect to $a$ and $x^\prime$($x^\prime$ is the
color appeared once in $c^\prime$ and $c^\prime(a)\neq x^\prime$. In
Cases 2 and 4 of Definition \ref{balanced}, we choose $x$ as the
color $x^\prime$. Then we have:

(a) The balanced coloring $c$ of $G:=H\bigcup P$ from $c^\prime$ is
a revised $\lceil\frac{|G|}{2}\rceil$-rainbow vertex-coloring such
that every color appears at most twice.

(b) When $|G|$ is odd, for any vertex $v\in V(G)$ there exists a balanced
coloring $c$ of $G$ satisfying the property ($*$) with respect to $v$ and $x_{\lceil
s/2\rceil-1}$($x_{\lceil s/2\rceil-1}$ is the color appeared once in $c$).
\end{pro}

\begin{proof}
We first prove (a). From the definition of balanced coloring and the
properties of $c^\prime$, the balanced coloring $c$ uses
$\lceil\frac{|G|}{2}\rceil$ colors and every color appears at most
twice. Since $c^\prime$ is a revised rainbow vertex-coloring of $H$,
$H$ is revised rainbow vertex-connected with respect to the balanced
coloring $c$. Hence, for any two vertices in $V(H)$ there exists a
revised rainbow path in $H$ with respect to $c$. Consider two
vertices $v_1\in V(H)$ and $v\in V(P)$. Let $P_a$ (resp. $P_b$) be a
revised rainbow path from $v_1$ to $a$ (resp. $b$) in $H$ such that
$x^\prime\notin colours(P_a)$ if $|H|$ is odd and $c(v_1)\neq
x^\prime$. Such a revised rainbow path $P_a$ exists since $c^\prime$
satisfies the property ($*$) with respect to $a$ and $x^\prime$ when
$|H|$ is odd. Then one of $P_a(aPv_2)$ and  $P_b(bPv_2)$ is a
revised rainbow path from $v_1$ to $v_2$. Note that when $|H|$ is
odd, the vertex on $P$ colored by $x^\prime$ has a revised rainbow
path to any vertex in $V(H)$ since $c^\prime$ satisfies the property
($*$) with respect to $a$ and $x^\prime$. Let $P^\prime$ be a
revised rainbow path from $a$ to $b$ in $H$. Then for any two
vertices in $V(P)$ there exits a revised rainbow path between them
on the cycle $P\bigcup P^\prime$. Therefore, the balanced coloring
$c$ of $G$ is a revised $\lceil\frac{|G|}{2}\rceil$-rainbow
vertex-coloring such that every color appears at most twice.

Now we prove (b). It is obvious that $|G|$ is odd in Cases 3 and 4
of Definition \ref{balanced}. First, consider $v\in
V(G)\setminus\{a, b\}$. In Case 3 of Definition \ref{balanced},
assign the color $x_{\lceil s/2\rceil-1}$ to the vertex $v_{\lceil
s/2\rceil}$. In Case 4 of Definition \ref{balanced}, if $c(v)\neq
c^\prime(a)$, then assign the color $x_{s/2-1}$ to $v_{s/2+1}$;
otherwise, we have $c(v)\neq c^\prime(b)$ and assign the color
$x_{s/2-1}$ to $v_{s/2}$. The other vertices are colored according
to Definition \ref{balanced}. From (a), the obtained balanced
coloring $c$ is a revised rainbow vertex-coloring. Let $P_a$(resp.
$P_b$) be a revised rainbow path from $v$ to $a$(resp. $b$) such
that $x^\prime\notin colours(P_a)$ if $|H|$ is odd and $c(v)\neq
x^\prime$. Then for any vertex $v^\prime\in V(H)$, the revised
rainbow path $P^\prime$ from $v$ to $v^\prime$ in $H$ satisfies that
$x_{\lceil s/2\rceil-1}\notin colours(P^\prime)$. For any vertex
$v^\prime\in V(P)$($c(v^\prime)\neq x_{\lceil s/2\rceil-1}$), one of
$P_a(aPv_2)$ and $P_b(bPv_2)$ is a revised rainbow path from $v$ to
$v^\prime$ such that $x_{\lceil s/2\rceil-1}$ does not appear on the
path. Second, consider $v\in V(P)$. Assign the color $x_{\lceil
s/2\rceil-1}$ to the vertex $u\in V(P)$ such that $d_P(u, v)=\lceil
s/2\rceil$ in Definition \ref{balanced}. It can be checked that the
obtained balanced coloring of $G$ has the property ($*$) with
respect to $v$ and $x_{\lceil s/2\rceil-1}$, i.e., for any vertex
$v^\prime\neq u$, there exists a revised rainbow path $P$ from $v$
to $v^\prime$ such that $u\notin V(P)$. It is clear that the color
$x_{\lceil s/2\rceil-1}$ appears once in $c$.
\end{proof}

Let $G$ be a 2-connected graph of order $n (n\geq3)$. Then $G$ has a
nonincreasing ear decomposition $(G_0, G_1, \cdots , G_k)$
satisfying the following conditions: (i) $G_0$ is an even cycle, if
$G$ is not an odd cycle; (ii) $G_i=G_{i-1}\bigcup P_i(1\leq i\leq
k)$ where $P_i$ is a longest ear of $G_{i-1}$, i.e., $\ell (P_1)\geq
\ell (P_2)\geq \cdots \geq \ell (P_k)$; (iii) $V(P_i)\bigcap
V(G_{i-1})=\{a_i, b_i\}(1\leq i\leq k)$ such that $a_i\neq b_i$. In
the sequel all the nonincreasing ear decomposition of a 2-connected
graph is one defined as above and the order of $G_i$ is denoted by
$n_i$. Without loss of generality, assume that $\ell (P_t)\geq 2$
and $\ell(P_{t+1})= \cdots =\ell(P_k)=1$. Hence, $G_t$ is a minimal
2-connected spanning subgraph of $G$, and every 2-connected graph
has such a spanning subgraph. It is obvious that $rvc(G_t)\leq
rvc(G)$. Therefore, in order to give an upper bound for the rainbow
vertex-connection of 2-connected graphs, we just need to consider
minimal 2-connected graphs.

\begin{lem} \label{lem1}
Let $G$ be a 2-connected graph of order $n(n\geq16)$. If the
nonincreasing ear decomposition ($G_0, G_1,\cdots, G_t$) of $G$
satisfies that $\ell(P_1)\geq\cdots\geq\ell(P_t)\geq 5$ (if
$t\geq1$), then $G$ has a revised $\lceil\frac n{2}\rceil$-rainbow
vertex-coloring such that every color appears at most twice, i.e.,
$rvc(G)\leq rvc^*(G)\leq\lceil\frac n{2}\rceil$.
\end{lem}

\begin{proof}
We prove the lemma by demonstrating a revised $\lceil\frac n
{2}\rceil$-rainbow vertex-coloring of $G$ such that every color
appears at most twice. If $G$ is an odd cycle, i.e., $G=C_n=v_1,
v_2, \cdots, v_n,$ $v_{n+1}(=v_1)$, then define a vertex-coloring of
$G$ by $c(v_i)=x_i$ for $1\leq i\leq \lceil n/2\rceil$ and
$c(v_i)=x_{i-\lceil n/2\rceil}$ if $\lceil n/2\rceil+1\leq i\leq n$.
Since for any two vertices $u, v$ of $G$ the path $P$ between $u$
and $v$ with length $d_G(u,v)$ is a revised rainbow path, $c$ is a
revised $\lceil\frac n{2}\rceil$-rainbow vertex-coloring such that
every color appears at most twice.

In the following, we assume that $G$ is not an odd cycle. We will
apply Proposition \ref{pro1} to show that $G_i(0\leq i\leq t)$ has a
revised $\lceil\frac{n_i}{2}\rceil$-rainbow coloring $c_i$ such that
every color appears at most twice. Assume that $G_0=C_{n_0}=v_1,
v_2, \cdots, v_{n_0}, v_{n_0+1}(=v_1)$. Define a vertex-coloring
$c_0$ of $G_0$ by $c_0(v_i)=x_i$ for $1\leq i\leq n_0/2$ and
$c_0(v_i)=x_{i-n_0/2}$ if $n_0/2+1\leq i\leq n_0$. It can be checked
that $c_0$ is a revised $\frac{n_0}{2}$-rainbow vertex-coloring such
that every color appears twice. If $t=0$, the result holds. Assume
that $t>0$. Since the vertex-coloring of $c_0$ satisfies the
conditions of the Proposition \ref{pro1}, $G_1=G_0\bigcup P_1$ has a
balanced coloring $c_1$ from $c_0$ satisfying the properties: (a)
$c_1$ is a revised $\lceil\frac{n_1}{2}\rceil$-rainbow
vertex-coloring such that every color appears at most twice; (b)
when $n_1$ is odd and $t>1$, $c_1$ has the property ($*$) with
respect to $a_2$ and $x_1^\prime$($x_1^\prime$ is the color appeared
once in $c_1$). If $t=1$, the result holds. Consider the case that
$t\geq 2$. It is obvious that the balanced coloring $c_1$ of $G_1$
satisfies the conditions of Proposition \ref{pro1}. Hence, using
Proposition \ref{pro1} $t$ times we can obtain a balanced coloring
$c_i$ of $G_i(1\leq i\leq t)$ from $c_{i-1}$ satisfying the
properties: (a) $c_i$ is a revised
$\lceil\frac{n_i}{2}\rceil$-rainbow vertex-coloring such that every
color appears at most twice; (b) when $n_i$ is odd and $i<t$, $c_i$
has the property ($*$) with respect to $a_{i+1}$ and $x_i^\prime$
($x_i^\prime$ is the color appeared once in $c_i$). Therefore, we
obtain a revised $\lceil\frac{n}{2}\rceil$-rainbow vertex-coloring
$c_t$ of $G$ such that every color appears at most twice.
\end{proof}

\begin{thm}
Let $G$ be a 2-connected graph of order $n (n\geq3)$. Then
$$
rvc(G)\leq\left\{
\begin{array}{ll}
0 & $if $ n=3; \\
1 & $if $ n=4, 5; \\
3 & $if $ n=9; \\
\lceil\frac{n}{2}\rceil-1 & $if $ n=6, 7, 8, 10, 11, 12, 13 $ or $ 15; \\
\lceil\frac{n}{2}\rceil & $if $ n\geq 16 $ or $ n=14,
\end{array}
\right.
$$
and the upper bound is tight, which is achieved by the cycle $C_n$.
\end{thm}

\begin{proof}
For $3\leq n\leq15$, it can be checked that $rvc(G)\leq rvc(C_n)$,
i.e., the result holds from Theorem \ref{thm:cycle}. In the
following, we show that $rvc(G)\leq \lceil\frac n {2}\rceil$ for
$n\geq16$. Without loss of generality, assume that $G$ is a minimal
2-connected graph. So the nonincreasing ear decomposition $(G_0,
G_1, \cdots ,G_k)$ of $G$ satisfies that
$\ell(P_1)\geq\cdots\geq\ell(P_k)\geq2$ if $k\geq1$. If $k=0$ or
$\ell(P_1)\geq\cdots\geq\ell(P_k)\geq5$, then $G$ has a revised
$\lceil\frac n{2}\rceil$-rainbow vertex-coloring from Lemma
\ref{lem1}. Hence, $rvc(G)\leq rvc^*(G)\leq\lceil\frac n{2}\rceil$.

Now assume that $k\geq1$, $5\leq\ell(P_t)\leq\cdots\ell(P_1)$ and
$2\leq\ell(P_k)\leq\cdots\leq\ell(P_{t+1})\leq4$ ($t< k$). From
Proposition \ref{lem1}, $G_t$ has a revised
$\lceil\frac{n_t}{2}\rceil$-rainbow vertex-coloring $c_t$ such that
every color appears at most twice. Let $x$ be a color of $c_t$ and
$x_j(t+1\leq j\leq k)$ and $x_0$ are new colors.

\begin{figure}[h,t,b,p]
\begin{center}
\scalebox{0.5}[0.5]{\includegraphics{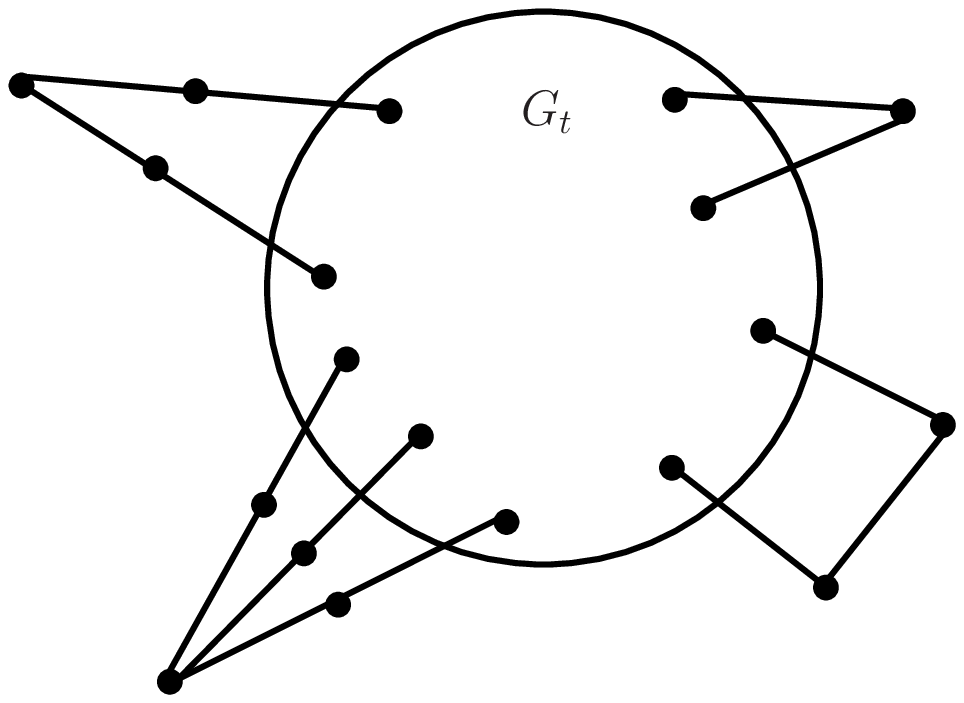}}

Figure 2. The graph used in the proof of Theorem 2.2.
\end{center}
\end{figure}

The graph $G$ is shown in Fig. 2 where the ears exist possibly.
Define a vertex-coloring $c$ of $G$ from $c_t$ as follows. For any
$v\in V(G_t)$, $c(v)=c_t(v)$. If there exists only one ear, say
$P_j=a_j, v_{j_1}, v_{j_2}, v_{j_3}, b_j (j=t+1)$ with length 4,
then $c(a_j)=c(v_{j_3})=x_j$, $c(v_{j_1})=c_t(a_j)$ and
$c(v_{j_2})=x$. If there exist at least two ears with length 4 and
$P_j=a_j, v_{j_1}, v_{j_2}, v_{j_3}, b_j (t+1\leq j\leq k)$ is such
an ear with length 4, then $c(a_j)=c(v_{j_3})=x_j$,
$c(v_{j_1})=c_t(a_j)$ and $c(v_{j_2})=x_0$. Note that the center
vertices of all ears with length 4 are colored by the new color
$x_0$. If $P_j= a_j, v_{j_1}, v_{j_2}, b_j(t+1\leq j\leq k)$ with
length 3, then $c(a_j)=c(v_{j_2})=x_j$ and $c(v_{j_1})=c_t(a_j)$. If
$P_j=a_j, v_{j_1}, b_j(t+1\leq j\leq k)$ with length 2, then
$c(v_{j_1})=x$. Note that in the definition of coloring $c$ of $G$
above, some vertex $a_j(t+1\leq j\leq k)$ is possibly colored more
than once. If possible, we choose the new color $x_j(t+1\leq j\leq
k)$ as its color. Hence, we obtain a vertex-coloring $c$ of $G$ from
$c_t$.

It is clear that $c$ uses at most $\lceil\frac n{2}\rceil$ colors.
In the following, we show that $G$ is rainbow vertex-connected.
Since $c_t$ is a revised rainbow coloring of $G_t$, $G_t$ is revised
rainbow vertex-connected respect to the coloring $c$. Hence, for any
two vertices in $V(G_t)$ there exists a revised rainbow path between
them in $G_t$. For any two vertices $v_1\in V(P_j)\backslash
V(G_t)(t+1\leq j\leq k)$ and $v_2\in V(G_t)$, one of
$(v_1P_ja_j)P^\prime$ and $(v_1P_jb_j)P^{\prime\prime}$ where
$P^\prime$ (resp. $P^{\prime\prime}$) is a revised rainbow path from
$a_j$ (resp. $b_j$) to $v_2$ in $G_t$ is a rainbow path from $v_1$
to $v_2$. Consider two vertices $v_1, v_2\in V(G)\backslash V(G_t)$.
If $d_G(v_1, v_2)\leq2$, then there is a rainbow path with length no
more than 2 from $v_1$ to $v_2$ trivially. If $d_G(v_1, v_2)\geq3$,
assume that $v_1\in V(P_{j_1})$ and $v_2\in V(P_{j_2})(t+1\leq
j_1<j_2\leq k)$. Since $\ell(P_{j_2})\leq4$, one of
$a_{j_2}P_{j_2}v_2$ and $b_{j_2}P_{j_2}v_2$ (say
$a_{j_2}P_{j_2}v_2$) has length no more than 2. Let
$P_{j_1}^\prime$(resp. $P_{j_1}^{\prime\prime}$) be a revised
rainbow path from $a_{j_1}$(resp. $b_{j_1}$) to $a_{j_2}$ in $G_t$.
Then one of $(v_1P_{j_1}a_{j_1})P_{j_1}^\prime(a_{j_2}P_{j_2}v_2)$
and $(v_1P_{j_1}b_{j_1})P_{j_1}^{\prime\prime}(a_{j_2}P_{j_2}v_2)$
is a rainbow path from $v_1$ to $v_2$. Therefore, $c$ is a rainbow
coloring of $G$, i.e., $rvc(G)\leq\lceil\frac n{2}\rceil$.

From Theorem \ref{thm:cycle}, the upper bound is tight.
\end{proof}

\begin{thm}
\label{thm:block} Let $G$ be a connected graph. If $G$ has a block
decomposition $B_1, B_2, \cdots, B_k$ and $t$ cut vertices, then
$rvc(G)\leq rvc(B_1)+rvc(B_2)+\cdots +rvc(B_k)+t$.
\end{thm}

\begin{proof}
If $G$ is a complete graph, then $rvc(G)=0$ and the result holds trivially. If
$G$ is not a complete graph, then $rvc(G)\geq 1$. We prove the result
by demonstrating a rainbow vertex-coloring using at most
$\Sigma_{i=1}^krvc(B_i)+t$ colors. If all the blocks are complete
graphs, then define a vertex-coloring $c$ of $G$ as follows.
The $t$ cut vertices are colored by $t$ colors $x_1, x_2, \cdots,
x_t$, and other vertices are colored by $x_1$. In this case, for any two
vertices there exists a path $P$ between them whose interval vertices are cut
vertices, i.e., $P$ is a rainbow path. Hence $c$ is a rainbow vertex-coloring.
The result holds.

Now assume that $B_1, \cdots, B_s(s\geq1)$ are not complete and $B_{s+1},
\cdots, B_k$ are complete. Then there exists a rainbow vertex-coloring
$c_i$ of $B_i(1\leq i\leq s)$ using $rvc(B_i)$ colors such that
$colours(B_i)\bigcap colours(B_j)=\emptyset (1\leq i< j\leq s)$.
Define a rainbow vertex-coloring $c_i$ of $B_i(s+1\leq i\leq k)$ that every
vertex of $B_i$ is colored by a color appeared in $B_1$. We define a
vertex-coloring $c$ of $G$ as follows. For any $v\in V(B_i)(1\leq i\leq k)$
which is not a cut vertex, $c(v)=c_i(v)$ and the $t$ cut vertices are
colored distinct with $t$ new colors. It is easy to check that the
vertex-coloring $c$ of $G$ using at most $\Sigma_{i=1}^krvc(B_i)+t$ colors
is a rainbow vertex-coloring. Therefore, $rvc(G)\leq rvc(B_1)+rvc(B_2)+\cdots +rvc(B_k)+t$.
\end{proof}

\end{document}